\documentclass{elsart3-1}

\usepackage{amssymb}
\usepackage{amsfonts}
\usepackage{amssymb}
\usepackage{color}
\usepackage{amsmath}

\usepackage[english,francais]{babel}

\newtheorem{theorem}{Theorem}[section]

\newtheorem{e-proposition}[theorem]{Proposition}

\newtheorem{e-definition}[theorem]{Definition\rm}

\newtheorem{theoreme}{Th\'eor\`eme}[section]

\newtheorem{proposition}[theoreme]{Proposition}

\newtheorem{definition}[theoreme]{D\'efinition\rm}

\renewcommand{\theequation}{\arabic{equation}}
\def\og{\leavevmode\raise.3ex\hbox{$\scriptscriptstyle\langle\!\langle$~}}
\def\fg{\leavevmode\raise.3ex\hbox{~$\!\scriptscriptstyle\,\rangle\!\rangle$}}
\def\dbN{{\mathop{\rm l\negthinspace N}}}
\def\dbR{{\mathop{\rm l\negthinspace R}}}
\def\ds{\displaystyle}

\journal{the Acad\'emie des sciences}
\begin{document}
% place in the next line the header (rubrique) chosen for your article,
% if you know it (you can also have 2, format : Header1/Header2
\centerline{}
\begin{frontmatter}

% Title, authors and addresses

% use the thanksref command within \title, \author or \address for footnotes;
% use the ead command for the email address,
% and the form \ead[url] for the home page:
% \title{Title\thanksref{label1}}
% \thanks[label1]{}
% \author{Name\thanksref{label2}}
% \ead{email address}
% \ead[url]{home page}
% \thanks[label2]{}
% \address{Address\thanksref{label3}}
% \thanks[label3]{}
\selectlanguage{english}
\title{Finite Codimensional  Controllability for Evolution Equations}

% use optional labels to link authors explicitly to addresses:
% \author[label1,label2]{}
% \address[label1]{}
% \address[label2]{}
% The [label1] can be suppressed if there is only one address for all authors

\selectlanguage{english}
\author[authorlabel1]{Xu  Liu},
\ead{liux216@nenu.edu.cn.}
\author[authorlabel2]{Qi L\"u},
\ead{lu@scu.edu.cn.}
\author[authorlabel2]{Xu Zhang}
\ead{zhang\_xu@scu.edu.cn.}

\address[authorlabel1]{Key Laboratory of Applied Statistics of MOE,  School of Mathematics and Statistics, Northeast Normal
University, Changchun 130024, China.}
\address[authorlabel2]{School of Mathematics,  Sichuan
University, Chengdu 610064,  China.}

% If you know the dates of reception, and acceptation you can put them now;
%  idem the name of the person presenting the Note

\medskip
\begin{center}
{\small Received *****; accepted after revision +++++\\
Presented by £££££}
\end{center}

\begin{abstract}
\selectlanguage{english}
% Text of abstract in English
Motivated by infinite-dimensional
optimal control problems with endpoint state
constraints, in this Note, we introduce the notion of finite  codimensional   exact
controllability for evolution equations. It is shown that this new controllability is equivalent to the
finite codimensionality condition in the literatures to guarantee Pontryagin's maximum
principle. As examples,  LQ  problems  with
fixed endpoint state constraints for a wave and a heat equation are analyzed,
respectively. {\it To cite this article: A.
Name1, A. Name2, C. R. Acad. Sci. Paris, Ser. I
340 (2005).}

\vskip 0.5\baselineskip

\selectlanguage{francais}
% Text of abstract in French
\noindent{\bf R\'esum\'e} \vskip
0.5\baselineskip \noindent {\bf Contr\^{o}le
exacte co-dimensionnel fini pour des
\'{e}quations d'\'{e}volution. } Motiv\'{e} par
des probl\`{e}mes de contr\^{o}le optimal en
dimension infinie avec des contraintes sur
l'\'{e}tat final, nous introduisons la notion de
contr\^{o}le exacte co-dimensionnel fini pour
des \'{e}quations d'\'{e}volution. On
d\'{e}montre que cette nouvelle notion de
contr\^{o}labilit\'{e} est \'{e}quivalente \`{a}
la condition de codimensionalit\'{e} finie qui
garantit que le principe maximal de Pontryagin
n'est pas trivial. A titre d'exemple, les
probl¨¨mes LQ avec des contraintes d'\'{e}tat de
point d'extr\'{e}mit\'{e} fixes sont
analys\'{e}s pour l'\'{e}quation des ondes et
l'\'{e}quation de chaleur respectivement. {\it
Pour citer cet article~: A. Name1, A. Name2, C.
R. Acad. Sci. Paris, Ser. I 340 (2005).}

\end{abstract}
\end{frontmatter}

% now the Version française abrégée, if it exists
\selectlanguage{francais}
\section*{Version fran\c{c}aise abr\'eg\'ee}
% Text of your Version française abrégée here.
% Note you do not need to repeat here equations that you use in the
% main text - for example 'voir (3)' is quite acceptable.

\selectlanguage{english}
% main text
\section{Introduction}

% etc, etc

% The Appendices part is started with the command \appendix;
% appendix sections are then done as normal sections
% \appendix

% \section{}
% \label{}

% The Acknowledgements are an un-numbered section
%\section*{Acknowledgements}
% Acknowledgements text here

Controllability is one of the fundamental issues in
control theory.    Up to now,  there are
numerous works devoted to  controllability
problems of linear and nonlinear distributed
parameter systems. In this Note, we will introduce a new concept  on the finite codimensional exact
controllability  for linear control
systems.

Let $Y$, $Z$ and $U$ be Hilbert spaces. Denote
by $\mathcal{L}(Z; Y)$ the set of all bounded
linear operators from $Z$ to $Y$, by $Y^*$ the
dual  space of $Y$, by span$D$ the closed
subspace spanned by a subset $D$ of $Y$, and by
$\overline{\mbox{co}}D$ the convex closed hull
of $D$. We identify $U^*$ with $U$. Let  $T>0$  and $p\in (1,\infty]$. Write $\mathcal{U}_p=L^p(0,
T; U)$. Consider the following
linear control  system:
\begin{eqnarray}\label{11}
&y_t(t)=Ay(t)+F(t)y(t)+B(t)u(t),
\quad t\in (0, T]\qquad\mbox{and}\qquad y(0)=y_0,
\end{eqnarray}
where  $u$ is the control variable and $y$ is
the state variable,    $A\!: \mathcal{D}(A)\!\subset\!
Y\!\rightarrow\!Y$ is a linear operator generating
a $C_0$-semigroup on $Y$,  $F(\cdot)\in L^\infty(0, T;
\mathcal{L}(Y; Y))$, $B(\cdot)\in L^\infty(0,
T; \mathcal{L}(U; Y))$, and $y_0\in Y$. For any $y_0\in Y$ and $u(\cdot)\in \mathcal{U}_p$,
(\ref{11}) admits a unique  mild solution $y(\cdot)\equiv y(\cdot; y_0, u(\cdot))\in
C([0, T]; Y)$. Define the reachable   set of
(\ref{11}) as follows:
\begin{equation}\label{0}
\mathcal{R}(T; y_0)=\big\{ y(T; y_0, u(\cdot))\in Y\
\big|\ y(\cdot) \mbox{ is the mild solution of  }
(\ref{11})\mbox{ with some } u(\cdot)\in \mathcal{U}_p
\big\}.
\end{equation}

Let us first recall the notion of finite
codimensionality.

\begin{definition}\label{10.25-def1}
A subset $M$ of $Y$ is said to
be finite codimensional  in $Y$,   if
there exists an
 $x_0\in\overline{\mbox{co}} M$  so that
span$\big(M-\{x_0\}\big)$  is a finite
codimensional subspace of $Y$, and $\overline{\mbox{co}}\big(M-\{x_0\}\big)$ has at
least  one interior  point in this subspace.
\end{definition}

Now, we introduce the following  new  notion  of  finite
codimensional  exact   controllability for (\ref{11}).

\begin{definition}\label{d1}
System \eqref{11} is said to be finite
codimensional  exactly  controllable at the time
$T$,  if  $\mathcal{R}(T;0)$ is  a finite
codimensional subspace of $Y$.
\end{definition}

Recall that  $(\ref{11})$ is called exactly   controllable at the time $T$, if  $\mathcal{R}(T;  0)=Y$. Hence,  the finite  codimensional  exact controllability defined above is clearly weaker than the exact controllability.
In general,   the  finite codimensional exact
controllability   cannot be reduced to the usual
exact controllability problem. Indeed, this can be done only for the very special case that $A+F(t)$ in $(\ref{11})$ has an invariant
subspace $Y_0$,  which is   finite codimensional in $Y$ and independent of $t\in  [0, T]$.

The  finite codimensional  exact
controllability is motivated by the study of
optimal control problems with endpoint state
constraints for infinite-dimensional systems. It is well known that, as a necessary
condition for optimal controls, Pontryagin's
maximum principle was established  for very general
finite-dimensional systems
(\cite{PC}), which is one of the
milestones in control theory. Nevertheless, very surprisingly, it
fails for infinite-dimensional systems without further assumptions
(\cite{Egorov}, see also \cite{LY1}).   This leads to that  for a long time,  the Pontryagin
maximum principle had been studied only for
evolution equations without terminal state
constraints.
Until 1980s, by assuming
a suitable finite codimensionality condition,
\cite{Fattorini,ly,LY1}
obtained the Pontryagin-type maximum principle for
optimal control problems with endpoint
constraints. However, it is usually quite difficult
to verify this condition directly. In this Note, we reduce the finite codimensionality
condition to a suitable  finite
codimensional exact controllability problem.
By the duality technique,  such  a controllability problem   is further
reduced to some {\it a priori}  estimate for its dual problem, which
maybe is easily verified, at least for some nontrivial example.

We refer to \cite{llz} for a detailed proof of the results in
this Note and other related results.

\section{Main result}
Let  $\widetilde{\mathcal{U}}$ be  a  bounded subset of the Banach  space $\mathcal{U}_p$ and   $\overline{\mbox{co}}\ \widetilde{\mathcal{U}}$ have at least one
 interior  point. In the sequel, we choose
$$
M=\big\{ y(T)\in Y\ \big|\ y(\cdot)\mbox{ is the solution of } (\ref{11}) \hbox{ with }y_0=0\mbox{ and some } u(\cdot)\in\widetilde{\mathcal{U}} \big\}.
$$
Also,  we consider the following homogenous linear
equation:
\begin{eqnarray}\label{215}
\phi_t(t)=-A^*\phi(t)-F(t)^*\phi(t),
\quad t\in (0, T]\qquad\mbox{and}\qquad\phi(T)=\phi_T,
\end{eqnarray}
where $\phi_T\in X^*$, and
$A^*$ and $F(t)^*$ are  respectively the dual operators
of $A$  and $F(t)$. Denote by  $C$ a generic positive constant, and by $p'$ the H\"older conjugate of $p$. The main result of this
Note is as follows.

\begin{theorem}\label{t2}
  The following assertions are
equivalent:

{\bf (1)}  The set $M$ is finite
codimensional in $Y$;

{\bf (2)}  The  equation
$(\ref{11})$ is finite codimensional exactly
controllable in $Y$;

{\bf (3)} There is a
finite codimensional subspace $\widetilde Y\subseteq
Y^*$ so that  any  solution $\phi$ of
$(\ref{215})$ satisfies
$$
|\phi_T|_{Y^*}\leq C
|B(\cdot)^*\phi|_{L^{p'}(0, T; U)},\quad\quad\forall\ \phi_T\in \widetilde Y;
$$

{\bf  (4)} There is a  compact
operator $G$ from  $Y^*$ to a Banach space
$X$ so  that   any  solution $\phi$ of
$(\ref{215})$ satisfies
$$
|\phi_T|_{Y^*}\leq C\big[|B(\cdot)^*\phi|_{L^{p'}(0, T; U)} +
|G\phi_T|_{X}\big],  \quad\quad\forall\ \phi_T\in Y^*.
$$
\end{theorem}

 Theorem \ref{t2}  can be applied    to
study  optimal control problems  with endpoint
constraints for  nonlinear distributed parameter
systems. For concrete problems, as we shall see in the next section, one may use the fourth assertion in Theorem \ref{t2} to check the finite codimensional exact
controllability of
$(\ref{11})$.

\section{Two examples}

This section is devoted to checking the finite codimensionality conditions  in some LQ problems (with fixed endpoint constraints) for  a wave and heat equations.  Let   $\Omega$ be a
bounded domain in $\dbR^N$ (for some $N\in\dbN$)  with a smooth boundary $\Gamma$, and
$\omega$ be a nonempty open subset  of
$\Omega$. Denote by $\chi_\omega$ the characteristic function of $\omega$.  Consider the following controlled  wave  and heat equations:
\begin{eqnarray}\label{8}\left\{\begin{array}{lll}
&y_{tt}-\Delta y+a(x, t)y=\chi_\omega u  &\mbox{ in
}Q\equiv \Omega\times(0,T),\\
&y=0 &\mbox{ on }\Sigma\equiv \Gamma\times(0,T),\\
&y(0)=y_0, \ y_t(0)=y_1  &\mbox{ in }\Omega,
\end{array}
\right.
\quad\quad\mbox{ and  } \quad\quad
 \left\{
\begin{array}{lll}
&y_{t}-\Delta y=\chi_\omega u  &\mbox{ in
}Q,\\
&y=0 &\mbox{ on }\Sigma,\\
&y(0)=y^0  &\mbox{ in }\Omega,
\end{array}
\right.
\end{eqnarray}
where  $u\in L^2(Q)$ is  the control variable. In \eqref{8},  $(y_0, y_1)\in H^1_0(\Omega)\times L^2(\Omega)$,  $a(\cdot)\in L^\infty(Q)$
and $y^0\in L^2(\Omega)$ are given. Also,
for  given targets  $(z_0, z_1)\in H^1_0(\Omega)\times L^2(\Omega)$  and  $z^0\in L^2(\Omega)$,  set
$$
\begin{array}{rl}
&\mathcal{U}_{ad}^1=\big\{ u\in  L^2(Q)\ \big|\  \mbox{The  solution }y\mbox{ of  the  wave equation in }(\ref{8}) \mbox{ satisfies  }(y(T), y_t(T))=(z_0, z_1) \big\},\\
&\mathcal{U}_{ad}^2=\big\{ u\in  L^2(Q)\ \big|\  \mbox{The  solution }y\mbox{ of the  heat equation in }(\ref{8}) \mbox{ satisfies  } y(T)=z^0 \big\},
\end{array}
$$
and
$$
J(u(\cdot))=\ds\frac{1}{2}\displaystyle\int_Q \big[\alpha(x, t)y^2(x,  t)+ \chi_\omega \beta(x, t)u^2(x, t)\big]dxdt,
$$
where  $\alpha,  \beta\in L^\infty(Q)$ are two given functions.
Assume that  $\overline{u}_i$  is  an  optimal control, i.e., it satisfies that
$
J(\overline{u}_i(\cdot))=\inf\big\{ J(u(\cdot))\;\big|\;u\in \mathcal{U}_{ad}^i\big\}$  ($i=1, 2$).
Write $\mathbb{B}_1=\big\{u\in L^2(Q)\;\big|\; |u|_{L^2(Q)}\leq 1\big\}$ and
$$\begin{array}{rl}
&M_1=\big\{ (y(T),  y_t(T))\;\big|\;  y\mbox{ solves the wave equation in } (\ref{8})\mbox{ with }(y_0,  y_1)=(0,  0)\mbox{ and some } u\in \mathbb{B}_1 \big\},\\
&M_2=\big\{ y(T)\;\big|\;  y\mbox{ solves the heat equation in }(\ref{8})\mbox{ with } y^0=0\mbox{ and some } u\in \mathbb{B}_1 \big\}.
\end{array}
$$

Similar  to the analysis in \cite{2}, if the sets $M_1$ and $M_2$ are finite codimensional accordingly in $H^1_0(\Omega)\times L^2(\Omega)$ and $L^2(\Omega)$,   then
one can obtain nontrivial necessary conditions
for the optimal controls $\overline{u}_i$  $(i=1,  2)$. To verify this finite codimensionality condition,  let us consider the following backward wave and heat
equations:
  \begin{eqnarray}\label{4}\left\{
\begin{array}{lll}
\psi_{tt}-\Delta \psi+a(x, t)\psi=0  \qquad\qquad\qquad\quad\mbox{ in
}Q,\\
\psi=0 \qquad\qquad\qquad\qquad\qquad\qquad\qquad\quad\mbox{ on }\Sigma,\\
\big(\psi(T), \psi_t(T)\big)=\big(\psi_1,\psi_2\big)\in L^2(\Omega)\times H^{-1}(\Omega),
\end{array}
\right.
\quad\quad\mbox{ and  }\quad\quad
\left\{
\begin{array}{lll}
\varphi_{t}+\Delta \varphi=0  &\mbox{ in
}Q,\\
\varphi=0 &\mbox{ on }\Sigma,\\
\varphi(T)=\varphi_T\in L^2(\Omega).
\end{array}
\right.
\end{eqnarray}

By the fourth assertion in  Theorem \ref{t2} and the known observability inequality  for the wave equation in \eqref{4} with $a(\cdot)\equiv0$ (\cite{BLR}), we obtain  the following positive  result for the wave equation (with a rather general $a(\cdot)$) in \eqref{8}.
\begin{proposition}\label{llz11}
For any $a(\cdot)\in L^\infty(Q)$, if
$(\Omega,\omega,T)$ fulfills the geometric
optics condition  (see \cite{BLR}), then $M_1$
is finite codimensional in $H^1_0(\Omega)\times
L^2(\Omega)$.
\end{proposition}

By Theorem $\ref{t2}$,     Proposition $\ref{llz11}$ implies  that   under the geometric optics condition,  the wave equation in $(\ref{8})$   is finite codimensional  exactly controllable.  Notice that  under the same condition,   the exact controllability of the wave equation with a general coefficient $a(\cdot)$ is  still  an  open  problem.

Finally, by the fourth assertion in  Theorem \ref{t2} again and the contradiction argument,  we have the following negative result for the heat equation.

 \begin{proposition}\label{llz12}
For any $\Omega$, $\omega$ and $T>0$, $M_2$ is not finite
codimensional in $
L^2(\Omega)$.
\end{proposition}

By Proposition \ref{llz12}, the finite codimensionality condition fails for LQ problems for heat equations with fixed endpoint constraints.

\section*{Acknowledgements}

This work is partially supported by the NSF of China
under grants  11471231, 11371084, 11221101 and 11231007,   by the Fundamental Research Funds for the Central
Universities under grants 2015SCU04A02  and 2412015BJ011,  by the Fok Ying Tong Education Foundation under grant   141001,
by  PCSIRT under  grant IRT$\_$15R53,  and by Grant MTM2014-52347 of
the MICINN, Spain.

\end{document}